\documentclass[12pt]{amsart}

\usepackage{a4wide}
\usepackage{graphicx,eepic}
\usepackage{amsmath,amsfonts,amssymb}
\usepackage{times}
\usepackage{xy}
\usepackage{bm}

\newcommand{\la}{\lambda}
\newcommand{\be}{\beta}
\newcommand{\de}{\delta}
\newcommand{\De}{\Delta}
\newcommand{\al}{\alpha}
\newcommand{\ga}{\gamma}
\newcommand{\e}{\varepsilon}
\newcommand{\Om}{\Omega}

\newcommand{\Si}{\Sigma}
\newcommand{\BR}{\mathbb{R}}

\newcommand{\BN}{\mathbb{N}}
\newcommand{\A}{\mathcal{A}}
\newcommand{\U}{\mathcal{U}}

\newcommand{\R}{\mathfrak{R}}

\renewcommand{\L}{\mathcal{L}}

\renewcommand{\S}{\mathcal S}

\newcommand{\C}{\mathcal C}

\newcommand{\F}{\mathcal F}

\newcommand{\N}{\mathcal N}

\newcommand{\p}{\bm p}
\newcommand{\x}{\bm x}

\newcommand{\w}{\bm w}

\newtheorem{lemma}{Lemma}[section]

\newtheorem{prop}[lemma]{Proposition}
\newtheorem{thm}[lemma]{Theorem}

\theoremstyle{definition}

\newtheorem{example}[lemma]{Example}

\theoremstyle{remark}
\newtheorem{rmk}[lemma]{Remark}

\numberwithin{equation}{section} \numberwithin{table}{section}

\title[Combinatorics of linear IFS with overlaps]
{Combinatorics of linear iterated  \\
function systems with overlaps}
\author{Nikita Sidorov}
\address{School of Mathematics, The University of Manchester, P.O. Box
88, Sackville Street, Manchester M60 1QD, United Kingdom. E-mail:
sidorov@manchester.ac.uk}
\date{\today}

\begin{document}

\begin{abstract}
Let $\bm p_0,\dots,\bm p_{m-1}$ be points in $\BR^d$, and let
$\{f_j\}_{j=0}^{m-1}$ be a one-parameter family of similitudes of
$\BR^d$:
$$
f_j(\bm x) = \la\bm x + (1-\la)\bm p_j,\ j=0,\dots,m-1,
$$ where
$\la\in(0,1)$ is our parameter. Then, as is well known, there exists
a unique self-similar attractor $S_\la$ satisfying
$S_\la=\bigcup_{j=0}^{m-1} f_j(S_\la)$. Each $\x\in S_\la$ has at
least one address
$(i_1,i_2,\dots)\in\prod_1^\infty\{0,1,\dots,m-1\}$, i.e., $\lim_n
f_{i_1}f_{i_2}\dots f_{i_n}({\bf 0})=\x$.

We show that for $\la$ sufficiently close to 1, each $\x\in
S_\la\setminus\{\p_0,\dots,\p_{m-1}\}$ has $2^{\aleph_0}$ different
addresses. If $\la$ is not too close to 1, then we can still have an
overlap, but there exist $\x$'s which have a unique address.
However, we prove that almost every $\x\in S_\la$ has $2^{\aleph_0}$
addresses, provided $S_\la$ contains no holes and at least one
proper overlap. We apply these results to the case of expansions
with deleted digits.

Furthermore, we give sharp sufficient conditions for the Open Set
Condition to fail and for the attractor to have no holes.

These results are generalisations of the corresponding
one-dimensional results, however most proofs are different.
\end{abstract}

\maketitle

\section{One-dimensional case: an overview}
\label{sec:1D}

The purpose of this paper is to generalise certain results
concerning a one-dimensional model considered in \cite{GS, Sid-AMM,
Sid2} so we first describe this model. Let $\la\in(1/2,1)$ be our
parameter. Consider a pair of similitudes of $I=[0,1]$:
\begin{align*}
f_0(x)&=\la x,\\
f_1(x)&=\la x+1-\la.
\end{align*}
They constitute an {\em iterated function system} (IFS). More
precisely, choose 0 as a starting point, and for any sequence
$(\e_1,\e_2,\dots)$ of 0s and 1s:
\[
x=\lim_{N\to+\infty}f_{\e_1}\dots f_{\e_N}(0),
\]
a {\em forward iteration} (the limit is independent of our choice of
the starting point). The set of all $x$'s that are representable in
such a form, is called the {\em attractor} of an IFS. As is well
known, in this case $\la>1/2$ yields the attractor $[0,1]$.

Unlike a general IFS, in this model any composition of $f_0$ and
$f_1$ can be given in a very simple form:
\[
f_{\e_1}\dots f_{\e_N}(0)=(\la^{-1}-1)\sum_{n=1}^N\e_n\la^n,
\]
whence
\[
x=\lim_{N\to\infty}\sum_{n=1}^N\e_k\la^n=
(\la^{-1}-1)\sum_{n=1}^\infty\e_n\la^n
\]
(sometimes called a $\la$-{\em expansion} of $x$). Since $\la>1/2$,
we have a proper {\it overlap}
\[
f_0(I)\cap f_1(I)=[1-\la,\la],
\]
so one might expect a typical $x\in(0,1)$ to have infinitely many
distinct $\la$-expansions ($=$ {\em addresses})
$(\e_1,\e_2,\dots)\in\prod_1^\infty\{0,1\}$ -- which indeed proves
to be the case.

More precisely, put
\[
\R_{\la}(x)=\left\{(\e_n)_1^\infty :
x=(\la^{-1}-1)\sum_{n=1}^\infty\e_n\la^{n}\right\}.
\]

The following important result regarding $\R_{\la}(x)$ has been
obtained by Erd\H os, Jo\'o and Komornik:

\begin{thm}\label{thm:ejk} \cite{EJK} If
$\la>g=\frac{\sqrt5-1}2=0.618\dots$, then $\R_\la(x)$ has the
cardinality of the continuum for each $x\in(0,1)$.
\end{thm}

Moreover, the golden ratio $g$ proves to be a sharp constant in the
previous theorem. Nonetheless, the following metric result holds for
$\la\in(1/2,g]$:

\begin{thm} \cite{Sid-AMM}\label{thm:sid} For any
$\la\in(1/2,g]$ the cardinality of $\R_{\la}(x)$ is the continuum
for Lebesgue-a.e. $x\in(0,1)$.
\end{thm}

Put
\[
\U_\la=\left\{x\in(0,1)\mid !\ (\e_n)_1^\infty:
x=(\la^{-1}-1)\sum_{n=1}^\infty\e_n\la^{n}\right\}
\]
(the {\it set of uniqueness}). By Theorem~\ref{thm:ejk},
$\U_\la=\emptyset$ if $\la>g$.

\begin{thm}\cite{GS}\label{thm:gs}
The set $\U_\la$ is:
\begin{itemize}
\item countable for $\la\in(\la_*,g)$;
\item uncountable of zero Hausdorff dimension if
$\la=\la_*$; and
\item a set of positive Hausdorff dimension for
$\la\in (1/2,\la_*)$.
\end{itemize}
Here $\la_*=0.559525\dots$ denote the (transcendental)
Komornik-Loreti constant introduced in \cite{KL}.
\end{thm}

The purpose of this paper is to generalise some of these results to
linear IFSs in higher dimensions.

\section{Multidimensional case: an analogue of
Theorem~\ref{thm:ejk}}

Let $\bm p_0,\dots,\bm p_{m-1}$ be distinct points in $\BR^d$, and
let $\{f_j\}_{j=0}^{m-1}$ be a one-parameter family of similitudes
of $\BR^d$:
\begin{equation}\label{fj}
f_j(\bm x) = \la\bm x + (1-\la)\bm p_j,\ j=0,\dots,m-1,
\end{equation}
where $\la\in(0,1)$ is our parameter\footnote{To simplify our
notation, we have decided to avoid notation like $f_j^{(\la)}$,
since there is never really any confusion regarding which $\la$ is
considered at a given moment.}.

Then, as is well known, there exists a unique self-similar {\em
attractor} $S_\la$ satisfying
\[
S_\la=\bigcup_{j=0}^{m-1} f_j(S_\la).
\]
Put $\A=\{0,\dots,m-1\}$. Similarly to the one-dimensional model,
every $\x\in S_\la$ has at least one {\em address}, i.e., a sequence
$(i_1,i_2,\dots)\in\A^\BN$ such that
\begin{align*}
\bm x&=\lim_{n\to+\infty} f_{i_1}\dots f_{i_n}(\bm x_0)\\
&=(\la^{-1}-1)\sum_{n=1}^\infty \la^n\bm p_{i_n},
\end{align*}
where $\bm x_0\in\BR^d$ is arbitrary. We assume the dimension of the
convex hull of $\{\p_0,\dots,\p_{m-1}\}$ to be equal to $d$.
(Otherwise we embed $S_\la$ into $\BR^{d'}$ with $d'<d$.)

Recall that an IFS is said to satisfy the Open Set Condition (OSC)
if there exists an open set $O\subset\BR^d$ such that
\[
O=\bigcup_{j=0}^{m-1} f_j(O),
\]
with the union being disjoint. Loosely speaking, the OSC means that
the images $f_j(\Om)$ do not intersect properly, where $\Om$ is the
convex hull of the $\p_j$. Virtually all famous IFS-generated
fractals (the Sierpi\'nski gasket, Sierpi\'nski carpet, von Koch
curve, etc.) originate from IFSs that satisfy the OSC.

We will be interested in IFSs which do {\em not} satisfy the OSC.
Here is a simple sufficient condition:

\begin{prop}\label{prop:osc}
If $\la>m^{-1/d}$, then the OSC is not satisfied.
\end{prop}
\begin{proof}Assume there exists an open set $O$ which satisfies
the definition. Since $f_j(O)\subset O$ and the images are disjoint,
we have $f_if_j(O)\cap f_{i'}f_{j'}(O)=\emptyset$ if
$(i,j)\neq(i',j')$, and by induction, $f_{i_1}\dots f_{i_n}(O)\cap
f_{j_1}\dots f_{j_n}(O)=\emptyset$ provided
$(i_1,\dots,i_n)\neq(j_1,\dots,j_n)$.

Since $\L_d (f_{i_1}\dots f_{i_n}(O))=\la^{dn}\L_d(O)$ (where $\L_d$
denotes the $d$-dimensional Lebesgue measure), and the number of
different words of length~$n$ is $m^n$, the pigeonhole principle
yields a contradiction with $\la>m^{-1/d}$.
\end{proof}

\begin{example}
Let $\bm p_0,\bm p_1, \bm p_2$ be three noncollinear points in
$\BR^2$ -- vertices of a triangle $\De$. Consider the IFS
\[
f_j(\bm x)=\la\bm x+(1-\la)\bm p_j,\quad j=0,1,2,
\]
and, following \cite{BMS}, we denote the attractor by $\S_\la$,
i.e.,
\[
\S_\la=\bigcup_{j=0}^2 f_j(\S_\la).
\]
Note that for $\la=1/2$ the set $\S_\la$ is the famous Sierpi\'nski
gasket. If $\la\le1/2$, then the IFS does satisfy the OSC, and for
$\la\ge2/3$ we have $\S_\la=\De$, i.e., $\S_\la$ contains no holes.
If $\la>1/2$, then we have a {\em proper overlap}, i.e.,
$f_i(\De)\cap f_j(\De)$ has a nonempty interior.
\end{example}

\begin{rmk}For the triangular case a more delicate argument allows
one to show that $\la>1/2$ implies the failure of the OSC (instead
of $\la>1/\sqrt3$ provided by Proposition~\ref{prop:osc}) -- see
\cite[Proposition~3.9]{BMS}.
\end{rmk}

Return to the general case. Put
$$
\Om=\mathrm{conv}(\bm p_0,\dots,\bm p_{m-1}).
$$
Clearly, $S_\la\subset\Om$. We give a universal sufficient condition
for $\Om$ to have no holes.

\begin{prop}\label{prop:dd}
If $\la\ge d/(d+1)$, then $S_\la=\Om$, i.e., our attractor has no
holes.
\end{prop}

\begin{proof}Let $\bm p_0,\dots,\bm p_{k-1}$ be the vertices of
$\Om$ (with $k\le m$) and let $\F_1,\dots,\F_h$ denote its
$(d-1)$-dimensional faces. Notice that $k\ge d+1$, and if $k=d+1$,
we have a simplex for which the claim is proved in \cite{BMS}.
Assume $k>d+1$.

For $\x\in\Om$ we denote its distance to $\F_i$ by $x_i$. Then
adding the volumes of the pyramids with the vertex $\x$ and the
bases $\F_i$ yields
\begin{equation}\label{eq:vol1}
\frac1d\sum_{i=1}^h \L_{d-1}(\F_i)\cdot x_i=\L_d(\Om).
\end{equation}
It suffices to show that $\Om=\bigcup_{j=0}^{k-1} f_j(\Om)$. Assume,
on the contrary, that there exists
$\x\in\Om\setminus\bigcup_{j=0}^{k-1} f_j(\Om)$. Since
$f_j(\Om)=\la\Om+(1-\la)\bm p_j$, we have
\begin{equation}\label{eq:vol2}
\frac{x_i}{\mbox{dist}\,(\bm p_j,\F_i)}<1-\la,\quad i=1,\dots,h.
\end{equation}
for all $j\in\{0,\dots,k-1\}$ such that $\bm p_j\notin\F_i$. Put
\[
\al_i=\max_{0\le j\le k-1} \mbox{dist}\,(\bm p_j,\F_i)\quad
i=1,\dots,h.
\]
Then by (\ref{eq:vol1}) and (\ref{eq:vol2}),
\[
\la<1-\frac{d\L_d(\Om)}{\sum_{i=1}^h \al_i\L_{d-1}(\F_i)}.
\]
To complete the proof, it suffices to show that
\begin{equation}\label{eq:helly}
\sum_{i=1}^h \al_i\L_{d-1}(\F_i)\le d(d+1)\L_d(\Om).
\end{equation}
Consider the family of affine copies of $\Om$ with the ratio
$d/(d+1)$, i.e.,
$$
\left\{\frac d{d+1}\,\Om+\frac1d\ \bm p_j\mid
j\in\{0,1,\dots,k-1\}\right\}.
$$
Let $0\le j_1<\dots<j_{d+1}\le k-1$ and let $\De_{j_1\dots
j_{d+1}}=\mbox{conv}(\bm p_{j_1},\dots,\bm p_{j_{d+1}})$, a
$d$-dimensional simplex. We have
\[
\bigcap_{r=1}^{d+1} \left(\frac d{d+1}\,\De_{j_1\dots
j_{d+1}}+\frac1d\ \bm p_{j_r}\right)\neq\emptyset
\]
(the intersection contains the centre of mass of $\De_{j_1\dots
j_{d+1}}$), whence
\[
\bigcap_{r=1}^{d+1} \left(\frac d{d+1}\,\Om+\frac1d\ \bm
p_{j_r}\right)\neq\emptyset
\]
for any $(d+1)$-tuple, as $\De_{j_1\dots j_{d+1}}\subset\Om$. Hence,
by Helly's theorem (see, e.g., \cite{Eck}), there exists
\[
\bm z\in \bigcap_{j=0}^{k-1}  \left(\frac d{d+1}\,\Om+\frac1d\ \bm
p_{j}\right).
\]
By our construction, the point $\bm z$ has the following property:
if $\bm y\in\partial\Om$, and $\bm z\in[\bm p_j,\bm y]$ for some
$j\in\{0,1,\dots,k-1\}$, then $|[\bm z,\bm
y]|\ge\frac1{d+1}|[\p_j,\bm y]|$. Therefore, $\mbox{dist}\,(\bm
p_j,\F_i)\le (d+1)\cdot\mbox{dist}(\bm z,\F_i)$ for all
$j\in\{0,\dots,k-1\}$, whence by definition, $\al_i\le
(d+1)\cdot\mbox{dist}(\bm z,\F_i)$. Consequently,
\[
\sum_{i=1}^h \al_i\L_{d-1}(\F_i)\le (d+1)\sum_{i=1}^h
\mbox{dist}(\bm z,\F_i)\cdot \L_{d-1}(\F_i),
\]
and to obtain (\ref{eq:helly}), it suffices to note that the volume
of $\Om$ equals the sum of the volumes of pyramids whose vertex is
$\bm z$, i.e.,
\[
\L_d(\Om)=\frac1d\sum_{i=1}^h \mbox{dist}(\bm z,\F_i)\cdot
\L_{d-1}(\F_i).
\]
\end{proof}

\begin{lemma}\label{lem:first}
Assume $\la$ is such that $S_\la=\Om$ and suppose $i,j\in\A$ are
such that $\Om_{ij}:=f_i(\Om)\cap f_j(\Om)$ has a nonempty interior.

Then each $\x\in\Om_{ij}$ has at least one address beginning with
$i$ and at least one address beginning with $j$.
\end{lemma}
\begin{proof}Since $\x\in f_i(\Om)$, there exists $\x'\in\Om$ such
that $f_i(\x')=\x$. By our assumption, $\x'\in S_\la$, whence
$\x'=\lim_n f_{i_2}\dots f_{i_n}(\x_0)$. Therefore, in view of the
continuity of $f_i$, we have $\x=\lim_n f_if_{i_2}\dots
f_{i_n}(\x_0)$. The same argument applies to $j$.
\end{proof}

\begin{rmk}The condition $S_\la=\Om$, generally speaking, cannot be
dropped; for instance, one can show that in the triangular case, if
$\la\in(0.65,2/3)$, then we have a hole whose image under one of the
maps lies in $\S_\la$ -- see \cite[Proposition~3.7]{BMS}.
\end{rmk}

\begin{thm}\label{thm:multi-ejk}
For each $\bm p_0,\dots,\bm p_{m-1}$ there exists $\la_0<1$ such
that for any $\la\in(\la_0,1)$,
\begin{enumerate}
\item There are no holes, i.e., $S_\la=\Om$;
\item Each point $\bm x\in\Om$, except when $\bm x$ is a vertex of
$\Om$, has $2^{\aleph_0}$ distinct addresses.
\end{enumerate}
\end{thm}
\begin{proof}Assume that $\la$ is large enough to ensure
$S_\la=\Om$. We will show that each $\bm x$ under consideration has
a continuum of addresses $(i_1,i_2,\dots)\in\A'$, where $\A'$ is the
set of indices $i$ such that $\p_i\in\partial\Om$.

The idea of the proof is to use the multivalued {\it inverse map}
$T_\la=\{f_0,\dots,f_{m-1}\}^{-1}$. More precisely, put
\[
\Om_i=f_i(\Om)\setminus\bigcup_{j\neq i} f_j(\Om).
\]
Clearly, if $\Om_i\ni \x\sim(i_1,i_2,\dots)$, then necessarily
$i_1=i$. Conversely, if a point $\x\notin \Om_i$ for any $i$, then
by Lemma~\ref{lem:first}, there is a choice for the first symbol of
its address.

Let $\x\in \Om_i\setminus\{\p_i\}$; note that shifting its address
$(i,i_2,i_3,\dots)$ yields $(i_2,i_3,\dots)$, which in $\Om$
corresponds to applying $f_i^{-1}$. Let $\la_0$ be such that
$f_i^{-1}(\Om_i)\cap \Om_j=\emptyset$ for all $i\neq j$ and all
$\la>\la_0$.

Thus, by $f_i^{-1}$ being expanding on $\Om_i\setminus\{\p_i\}$, we
conclude that there exists $k\ge1$ such that
\begin{align*}
f_i^{-k}(\x) & \in \Om_i,\\
\x'=f_i^{-k-1}(\x) & \in \Om\setminus\bigcup_{i=0}^{m-1} \Om_i.
\end{align*}
By our construction, $\x'$ has at least two addresses; consider its
two shifts, $\x''$ and $\x'''$, say. Either
$\x''\in\Om\setminus\bigcup_{i=0}^{m-1} \Om_i$ and thus, has at
least two addresses itself or it belongs to $\Om_j$ for some $j$
(and obviously, is not equal to $\p_j$) and, similarly to the above,
we shift its address until it falls into
$\Om\setminus\bigcup_{i=0}^{m-1} \Om_i$. Hence any $\x\in\Om$ that
is not one of its vertices, has $2^{\aleph_0}$ distinct addresses.
\end{proof}

\begin{rmk}
Note that for the triangular case the sharp constant is
$\la_0\approx0.68233$, the unique positive root of $x+x^3=1$ -- see
Theorem~\ref{thm:triangle}.
\end{rmk}

\section{Multidimensional case: generic
behaviour}\label{sec:generic}

\subsection{General theory: Lebesgue measure.}
Let $U_\la$ denote the set of $\x\in S_\la$ having a
unique address, and $R_\la(\x)$ denote the set of all addresses of a
given $\x\in S_\la$.

\begin{lemma}\label{lem:less}
Put
\[
V_\la = \left\{\x\in S_\la :
\mathrm{card}\,R_\la(\x)<2^{\aleph_0}\right\}.
\]
Then $\dim_H V_\la=\dim_H U_\la$ (where $\dim_H$ denotes Hausdorff
dimension). 
\end{lemma}

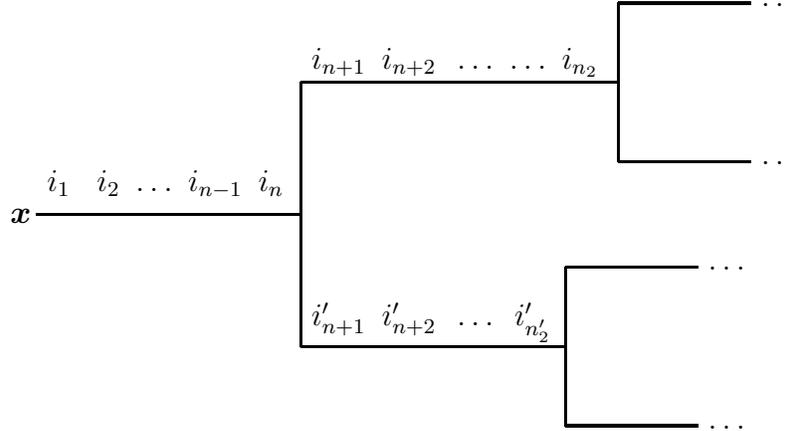
\begin{figure}[t]
    \centering 
\centerline{
\begin{picture}(265,210)
\thicklines  \put(0,100){\line(1,0){100}}
\put(100,100){\line(0,-1){50}}
\put(100,100){\line(0,1){50}}\put(100,150){\line(1,0){120}}
\put(100,50){\line(1,0){100}}
\put(220,150){\line(0,1){30}} \put(220,150){\line(0,-1){30}}
\put(200,50){\line(0,1){30}} \put(200,50){\line(0,-1){30}}
\put(200,20){\line(1,0){50}}\put(200,80){\line(1,0){50}}
\put(220,120){\line(1,0){50}}\put(220,180){\line(1,0){50}}
\put(274,179){$\dots$}\put(274,119){$\dots$}
\put(254,79){$\dots$}\put(254,19){$\dots$}
\put(4,109){$i_1\,\,\,\,\, i_2\,\,\dots\,\, i_{n-1}\,\,\, i_n$}
\put(104,155){$i_{n+1}\,\,\,i_{n+2}\,\,\,\dots\,\,\dots\,\,
i_{n_2}$}
\put(104,58){$i'_{n+1}\,\,\,i'_{n+2}\,\,\,\dots\,\,\,i'_{n_2'}$}
\put(-10,97){$\x$}
\end{picture}}
\caption{Branching and bifurcations.}
    \label{fig:branching}
\end{figure}

\begin{proof}We are going to exploit the idea of {\it branching}
introduced in \cite{Sid2}. Let $\x\in S_\la$ have at least two
addresses; then there exists the smallest $n\ge0$ such that
$\x\sim(i_1,\dots,i_n,i_{n+1},\dots)$ and
$\x\sim(i_1,\dots,i_n,i'_{n+1},\dots)$ with $i_{n+1}\neq i'_{n+1}$.
We may depict this {\em bifurcation} as is shown in
Fig.~\ref{fig:branching}.

Assume that $\x$ has less than a continuum of distinct addresses.
Then, inevitably, one of the branches at some point ceases to
bifurcate. In other words, there exists $(i_n)_1^\infty\in\A^\BN$
and $N\in\BN$ such that $\x\sim(i_1,i_2,\dots)$ and
$\x'\sim(i_N,i_{N+1},\dots)\in U_\la$. Hence
\begin{equation}\label{eq:Vla}
V_\la \subset \bigcup_{(i_1,\dots,i_N)\in\A^N} f_{i_1}\dots
f_{i_N}(U_\la).
\end{equation}
Since the $f_i$ are linear, (\ref{eq:Vla}) implies $\dim_H V_\la\le
\dim_H U_\la$, and the inverse inequality is trivial.
\end{proof}

\begin{rmk}
Note that if $\la_0$ in Theorem~\ref{thm:multi-ejk} is sharp, we
always have a nonempty set of uniqueness for $\la<\la_0$, because,
as we know from the branching argument, the existence of $\x$ with
less than a continuum of addresses implies the existence of $\x'$
with a unique address. For an example see
Section~\ref{sec:triangle}.
\end{rmk}

Our goal is to show that, similarly to the one-dimensional case, if
there are no holes and at least one proper overlap, then a.e. $\x$
has a continuum of addresses. We need an auxiliary claim from
dimension theory:

\begin{lemma}\label{lem:mesh}
Let $A\subset\BR^d$ be such that there exists a positive constant
$\de>0$ such that for an arbitrary cube $\C \subset \BR^d$ which
intersects $A$, one can find a cube $\C_0\subset C$ such that
$\L_d(\C_0)\ge\de\L_d(\C)$ and $\C_0\cap A=\emptyset$.

Then $\dim_HA<d$.
\end{lemma}
\begin{proof}Since $\dim_H A\le\underline{\dim}_B A$ (lower
box-counting dimension), it suffices to show that
$\underline{\dim}_B A<d$. Recall that there are various definitions
of $\underline{\dim}_B$ (see \cite[Chapter~3.1]{FG}) and in
particular, the one which involves mesh cubes which we will use.

More precisely, a cube of the form
$[m_1\e,(m_1+1)\e]\times\dots\times [m_d\e,(m_d+1)\e]$ for some
$\e>0$ and $m_1,\dots,m_d$ integers, is called an {\it $\e$-mesh
cube}. Let $\N_\e(A)$ denote the number of $\e$-mesh cubes which
intersect $A$. Then
\[
\underline{\dim}_B A=\liminf_{\e\to0}\frac{\log\N_\e(A)}{\log
(1/\e)}.
\]
It is obvious that our condition implies that there exists $M\in\BN$
such that for any $\e$-mesh cube $C$ which intersects $A$, there
exists an $\e/M$-mesh cube $C_0\subset C$ which doesn't. Hence
\[
\N_{\e/M}(A)\le (M^d-1)\N_\e(A),
\]
which implies
\[
\N_{M^{-n}}(A)\le (M^d-1)^n\cdot\N_1(A).
\]
Consequently,
\[
\liminf_{\e\to0}\frac{\log\N_\e(A)}{\log
1/\e}\le\frac{\log(M^d-1)}{\log M}<d.
\]
\end{proof}

Now we are ready to prove a key technical lemma.

\begin{lemma}\label{lem:mes0}
Let $\{F_j\}_{i=0}^{L-1}$ be a finite family of linear contractions
of $\BR^d$ with the same contraction ratio $\be\in(0,1)$. Let
$\Om\subset\BR^d$ be a polyhedron of dimension~$d$, and assume
$F_j(\Om)\subset\Om$ for all $j=0,\dots,L-1$, and
$\Om=\bigcup_{j=0}^{L-1} F_j(\Om)$ (no holes).

Put
\begin{equation}
\label{eq:snwn}
W_n=\bigcup_{\substack{(j_1,\dots,j_n)\in\{0,\dots,L-1\}^n:\\
\exists\ k\in\{1,\dots,n\}: j_k=0}}F_{j_1}\dots F_{j_n}(\Om).
\end{equation}
Then $W_n\subset W_{n+1}$ for all $n\ge1$, and $\L_d(\Om\setminus
W)=0$, where $W=\bigcup_{n\ge1} W_n$. Furthermore,
$\dim_H(\Om\setminus W)<d$.
\end{lemma}
\begin{proof}Our first goal is to show that $W_n\subset W_{n+1}$ for
all $n\ge1$. Since $\Om=\bigcup_{j=0}^{L-1} F_j(\Om)$, we have by
induction,
\[
\bigcup_{(j_1,\dots,j_k)\in\{0,\dots,L-1\}^k}F_{j_1}\dots
F_{j_k}(\Om)=\Om,
\]
whence for any $k<n$,
\[
\bigcup_{(j_1,\dots,j_n)\in\{0,\dots,L-1\}^n}F_{j_1}\dots F_{j_k}
F_0 F_{j_{k+1}}\dots
F_{j_n}(\Om)=\bigcup_{(j_{1},\dots,j_k)\in\{0,\dots,L-1\}^{k}}
F_{j_{1}}\dots F_{j_k}F_0(\Om).
\]
Hence by definition, $W_1=F_0(\Om)$ and
\[
W_n=\left(\bigcup_{k=1}^{n-1}
\bigcup_{(j_{1},\dots,j_k)\in\{0,\dots,L-1\}^{k}} F_{j_{1}}\dots
F_{j_k}F_0(\Om)\right)\cup F_0(\Om),\quad n\ge2.
\]
Consequently,
\[
W_{n+1}=W_n\cup\bigcup_{(j_1,\dots,j_n)\in\{0,\dots,L-1\}^n}
F_{j_1}\dots F_{j_n}F_0(\Om),\quad n\ge1,
\]
whence $W_n\subset W_{n+1}$.

By Lemma~\ref{lem:mesh}, to show that $\dim_H(\Om\setminus W)<d$, it
suffices to demonstrate that there exists a positive constant
$\de=\de(\Om,\be)>0$ such that given an arbitrary cube $\C \subset
\Om$, one can find a cube $\C_0\subset C$ such that
$\L_d(\C_0)\ge\de\L_d(\C)$ and $\C_0\cap (\Om\setminus
W)=\emptyset$.

So we choose an arbitrary cube $\C\subset\Om$ and denote the length
of its edge by $\kappa$. Let $\xi$ denote the centre of $\C$; since
$\Om$ has no holes, $\xi=\lim_{r\to\infty} F_{j_1}\dots
F_{j_r}(\Om)$ (the limit in the Hausdorff metric). Hence there
exists a unique $N$ such that $F_{j_1}\dots
F_{j_{N-1}}(\Om)\not\subset\C$, and $F_{j_1}\dots
F_{j_{N}}(\Om)\subset\C$.

Notice that since $\xi\in F_{j_1}\dots F_{j_{N-1}}(\Om)$ and
$F_{j_1}\dots F_{j_{N-1}}(\Om)\not\subset \C$, we have
$$
\mbox{diam}\,F_{j_1}\dots F_{j_{N-1}}(\Om)\ge\mbox{dist}\
(\xi,\partial \C)=\kappa/2.
$$
Put
\[
\nu=\frac{\mbox{diam}(\Om)}{(\L_d(\Om))^{1/d}}.
\]
Then
\begin{align*}
\L_d(\C)^{1/d}&=\kappa \le 2\,\mbox{diam}\,F_{j_1}\dots
F_{j_{N-1}}(\Om)=2\be^{N-1}\mbox{diam}(\Om)\\
&=\frac{2\nu}{\be}\cdot\be^N(\L_d(\Om))^{1/d}=\frac{2\nu}{\be}\cdot
(\L_d(F_{j_1}\dots F_{j_{N-1}}(\Om))^{1/d},
\end{align*}
whence
\begin{equation}\label{eq:ld}
\frac{\L_d(F_{j_1}\dots F_{j_{N}}(\Om))}{\L_d(\C)}\ge c>0,
\end{equation}
where $c=(\be/2\nu)^d$, i.e., $c$ depends only on the shape of $\Om$
and on the contraction ratio, but not on $\C$ itself.

Put $\Om_0=F_{j_1}\dots F_{j_N}F_0(\Om)$; by (\ref{eq:ld}),
$\L_d(\Om_0)\ge \be c\cdot \L_d(\C)$. We can find a cube
$\C_1\subset\Om$ such that the ratio of their volumes equals
$\ga>0$. Since $\Om_0$ is similar to $\Om$, we put
$\C_0=F_{j_1}\dots F_{j_N}F_0(\C_1)\subset\Om_0$ and obtain
\[
\L_d(\C_0)\ge \be\ga c\cdot\L_d(\C),
\]
where $\de:=\be\ga c$ is independent of $\C$. Furthermore, $\Om_0$
(and consequently, $\C_0$) has an empty intersection with
$\Om\setminus W_{N+1}$, whence $\C_0\cap(\Om\setminus W)=\emptyset$
as well, and we are done.
\end{proof}

\begin{thm}\label{thm:main}
Assume
\begin{itemize}
\item $S_\la=\Om$, i.e., there are no holes;
\item there exist $i,k\in\A$ such that a vertex of
$f_k(\Om)$ belongs to the interior of $f_i(\Om)$.
\end{itemize}
Then $\L_d$-a.e. $\x\in\Om$ has $2^{\aleph_0}$ distinct addresses,
and the exceptional set $V_\la$ has Hausdorff dimension strictly
less than~$d$.
\end{thm}
\begin{rmk}If $d\le2$, it suffices to assume that $f_i(\Om)\cap
f_k(\Om)$ has a nonempty interior, since if two convex polygons (or
intervals) intersect properly, then it is obvious that there exists
a vertex of one which lies in the interior of the other. For $d\ge3$
this is not always the case.
\end{rmk}
\begin{proof}By our assumption, there exists $j\in\A$ such
that $f_k(\p_j)\in int(f_i(\Om))$. Hence there exists $\ell\in\BN$
such that $f_if_j^{\ell-1}(\Om)\subset f_i(\Om)\cap f_k(\Om)$. Since
$S_\la=\Om$, Lemma~\ref{lem:first} implies that any $\x\in
f_if_j^{\ell-1}(\Om)$ has at least two different addresses.

Put $L=m^\ell$ and define $\{F_0,\dots,F_{L-1}\}=\{f_{i_1}\dots
f_{i_\ell}\mid (i_1,\dots,i_\ell)\in\A^\ell\}$ with
$F_0=f_if_j^{\ell-1}$. By the above, each $\x\in F_{j_1}\dots
F_{j_{k-1}}F_0 F_{j_{k+1}}\dots F_n(\Om)$ has at least two different
addresses, whence $U_\la\subset \Om\setminus W$, where $W=\bigcup_n
W_n$ and $W_n$ is given by (\ref{eq:snwn}).

Hence by Lemma~\ref{lem:mes0}, $\dim_H(U_\la)<d$, whence by
Lemma~\ref{lem:less}, $\dim_H(V_\la)<d$, which is the claim of the
theorem.
\end{proof}

\begin{rmk}If $\la$ is sufficiently close to the critical value
$\la_0$ (see the previous section), then one could expect the
exceptional set $V_\la$ to be countable (similarly to the
one-dimensional case).
\end{rmk}

\subsection{Application: $\la$-expansions with deleted digits.}
Expansions of real numbers in non-integer bases with deleted digits
have been studied since the mid-1990s -- see, e.g., \cite{KSS, PS}.
The model is as follows: assume $d=1$ and let
$A=\{a_1,\dots,a_m\}\subset\BR$ be a ``digit'' set with $a_1<\dots
<a_m$. Let $x\in\BR$ have an expansion of the form
\begin{equation}\label{eq:deleted}
x=\sum_{n=1}^\infty \e_n\la^n,\quad \e_n\in A,\ n\ge1.
\end{equation}
It is obvious that $\la a_1/(1-\la)\le x\le \la a_m/(1-\la)$.
M.~Pedicini \cite{Ped} has shown that if
\begin{equation}\label{eq:ped}
\max_{1\le j\le m-1} (a_{j+1}-a_j)<\frac{\la(a_m-a_1)}{1-\la},
\end{equation}
then each $x\in[\la a_1/(1-\la),\la a_m/(1-\la)]$ has at least one
expansion of the form (\ref{eq:deleted}). Note also that in the
recent paper \cite{DK} the theory of random and greedy
beta-expansions with deleted digits (under the
assumption~(\ref{eq:ped})) has been developed.

We apply our results from this and the previous section to obtain

\begin{prop}\label{prop:deleted}
\begin{enumerate}
\item There exists $\la_0=\la_0(a_1,\dots,a_m)<1$ such that for each
$\la\in(\la_0,1)$ any $x\in(\la a_1/(1-\la),\la a_m/(1-\la))$ has
$2^{\aleph_0}$ expansions of the form~(\ref{eq:deleted}).
\item If the condition~(\ref{eq:ped}) is satisfied, Lebesgue-a.e.
$x\in(\la a_1/(1-\la),\la a_m/(1-\la))$ has $2^{\aleph_0}$
expansions of the form~(\ref{eq:deleted}), and the exceptional set
has Hausdorff dimension strictly less than~1.
\end{enumerate}
\end{prop}
\begin{proof}Put
\begin{equation}\label{eq:ifs-del}
f_j(x)=\la (x+a_j),\quad 1\le j\le m.
\end{equation}
Then, as in the standard one-dimensional case (where $a_1=0,a_2=1$),
we have by induction,
\[
f_{\e_1}\dots f_{\e_n}(x_0)=\la^nx_0+\sum_{k=1}^n \e_k\la^k,
\]
whence
\[
\lim_{n\to\infty}f_{\e_1}\dots f_{\e_n}(x_0)=\sum_{n=1}^\infty
\e_n\la^n
\]
for any $x_0\in\BR$. Therefore, $x$ has an expansion of the
form~(\ref{eq:deleted}) if and only if $x\in S_\la$ for the
IFS~(\ref{eq:ifs-del}).

The condition~(\ref{eq:ped}) ensures that $S_\la=\Om=[\la
a_1/(1-\la),\la a_m/(1-\la)]$. To prove the first part of the
proposition, notice that (\ref{eq:ped}) holds for all $\la$
sufficiently close to 1 so Theorem~\ref{thm:multi-ejk} is applicable
to any $x$ which lies in the interior of $\Om$. (As
$\partial\Om=\{a_1/(1-\la),a_m/(1-\la)\}$.)

To prove the second part, we notice that by (\ref{eq:ped}),
\[
|f_j(\Om)|=\frac{\la^2(a_m-a_1)}{1-\la}>\la(a_{j+1}-a_j),\quad 1\le
j\le m-1.
\]
Hence
\begin{align*}
\sum_{j=1}^m |f_j(\Om)| &= \sum_{j=1}^{m-1} |f_j(\Om)|+|f_m(\Om)|\\
&> \la(a_m-a_1)+\frac{\la^2(a_m-a_1)}{1-\la}\\
&=\frac{\la(a_m-a_1)}{1-\la}=|\Om|,
\end{align*}
whence there exists $j\in\{1,\dots,m-1\}$ such that $f_j(\Om)\cap
f_{j+1}(\Om)$ has a nonempty interior. Thus, we can apply
Theorem~\ref{thm:main} to this setting.
\end{proof}

\begin{rmk}In her PhD dissertation, Anna-Chiara Lai \cite{Lai} has
proved a weaker version of the second claim of
Proposition~\ref{prop:deleted}.
\end{rmk}

Finally, we prove

\begin{lemma}Provided (\ref{eq:ped}) is satisfied, the Open Set
Condition for the IFS~(\ref{eq:ifs-del}) fails.
\end{lemma}
\begin{proof}We have
\[
a_m-a_1=\sum_{j=1}^{m-1}(a_{j+1}-a_j)<\frac{\la(m-1)(a_m-a_1)}{1-\la},
\]
whence $\la>1/m$, and we apply Proposition~\ref{prop:osc}.
\end{proof}

\subsection{General theory: natural measure.}
In the end of this section we would like to obtain a result similar
to Theorem~\ref{thm:main} for a ``natural'' measure on $S_\la$. Let
$(p_0,p_1,\dots,p_{m-1})$ be a probability vector with $p_j>0$ for
all $j$. The {\em probabilistic IFS} given by the $f_i$ and the
$p_i$ is defined as follows: put $\Si=\A^\BN$ and define $\rho$ as
the product measure on $\Si$ with equal multipliers
$(p_0,p_1,\dots,p_{m-1})$. Let the projection map $\pi:\Si\to\BR^d$
be given by the formula
\[
\pi(i_1,i_2,\dots):=\lim_{n\to+\infty}f_{i_1}f_{i_2}\dots
f_{i_n}(\mathbf 0).
\]
We define the measure $\mu$ on $S_\la$ as the push down measure
$\pi(\rho)$. As is well known, $\mbox{supp}(\mu)=S_\la$.

\begin{prop}\label{prop:pushdown}
Under the assumptions of Theorem~\ref{thm:main}, $\mu$-a.e.
$\x\in\Om$ has $2^{\aleph_0}$ distinct addresses.
\end{prop}
\begin{proof}Put $\w=ij^{\ell-1}$. By the Birkhoff ergodic theorem
applied to the one-sided Bernoulli shift on the measure space
$(\Si,\rho)$,
$$
\rho\{(i_1,i_2,\dots)\in\Si \mid \exists k :
(i_k,\dots,i_{k+\ell-1})=\w\}=1,
$$
whence $\mu(U_\la)=0$, because $\x\in U_\la$ cannot have an address
containing $\w$. All that is left is to show that $\mu(V_\la)=0$ as
well.

In view of (\ref{eq:Vla}), it suffices to show that
\begin{equation}\label{eq:mu}
\mu(f_{i_1}\dots f_{i_n}(U_\la))=0,\quad \forall
(i_1,\dots,i_n)\in\A^n.
\end{equation}
Note that, as is well known (see, e.g., \cite{DF}), the
self-similarity of the measure $\mu$ implies
\[
\mu(E)=\sum_{s=0}^{m-1} p_s\cdot \mu (f_s(E)),
\]
for any Borel set $E$. Hence by induction,
\[
\mu(E)=\sum_{(i_1,\dots,i_n)\in\A^n}p_{i_1}\dots p_{i_n}\cdot
\mu(f_{i_1}\dots f_{i_n}(E)),
\]
which implies (\ref{eq:mu}).
\end{proof}
\begin{rmk}In fact, to apply the ergodic theorem to the shift
on $\Si$, all we need from $\rho$ is $p_i>0$ and $p_j>0$; the other
components of the probability vector may equal zero.
\end{rmk}

The main problem for the future study is to check whether in some
cases of IFSs with holes an analogue of Theorem~\ref{thm:main} still
holds. We plan to be study this question in our subsequent papers.

\section{Main example: triangle}\label{sec:triangle}

For the triangular case we give an explicit analogue of
one-dimensional results mentioned in Section~\ref{sec:1D}. Following
\cite{BMS}, we denote the set of uniqueness by $\U_\la$.

\begin{thm}\label{thm:triangle}
Let $\la_0\approx0.68233$ be the unique positive root of $x^3+x=1$.
Then
\begin{enumerate}
\item for $\la<\la_0$, then the set of uniqueness $\U_\la$ is nonempty;
\item if $\la\in(\la_0,1)$, then each $x\in\S_\la
\setminus\{\p_0,\p_1,\p_2\}$ has $2^{\aleph_0}$ different addresses.
\end{enumerate}
\end{thm}
\begin{proof}(1) Firstly, we introduce a convenient coordinate system
for this case suggested in \cite{BMS}. Without loss of generality,
we may assume our triangle $\De$ to be equilateral. We now identify
each point $\x\in\De$ with a triple $(x,y,z)$, where
\[
x=\mbox{dist}\,(\bm x,[\p_1,\p_2]),\ y=\mbox{dist}\,(\bm
x,[\p_0,\p_2]),\ z=\mbox{dist}\,(\bm x,[\p_0,\p_1]),
\]
where $[\p_i,\p_j]$ is the edge containing $\p_i$ and $\p_j$.  As is
well known, $x+y+z$ equals the tripled radius of the inscribed
circle, and we choose it to be equal to 1. These coordinates are
called {\em barycentric}. Henceforward we write each $\x\in\S_\la$
in barycentric coordinates.

It is shown in \cite{BMS} that $(x,y,z)\in\S_\la$ if and only if
there exist three 0-1 sequences $(a_n)_0^\infty, (b_n)_0^\infty$ and
$(c_n)_0^\infty$ such that
\begin{align*}
x&=(1-\la)\sum_{n=0}^\infty a_n\la^n,\\
y&=(1-\la)\sum_{n=0}^\infty b_n\la^n,\\
x&=(1-\la)\sum_{n=0}^\infty c_n\la^n,
\end{align*}
with $a_n+b_n+c_n=1$ for all $n\ge0$.

We claim that the point
\begin{equation}\label{eq:pi}
\bm\pi(\la)=\left(\frac{\la^2}{1+\la+\la^2},
\frac{\la}{1+\la+\la^2}, \frac{1}{1+\la+\la^2}\right)
\end{equation}
belongs to $\U_\la$ provided $\la<\la_0$. To prove this, it suffices
to demonstrate that the system of equations
\begin{equation}
\begin{aligned}\label{eq:triple}
a_0+a_1\la+a_2\la^2+\dots&=\frac{\la^2}{1-\la^3},\\
b_0+b_1\la+b_2\la^2+\dots&=\frac{\la}{1-\la^3},\\
c_0+c_1\la+c_2\la^2+\dots&=\frac{1}{1-\la^3}
\end{aligned}
\end{equation}
has a unique solution $(a_n)_0^\infty=(001001\dots),
(b_n)_0^\infty=(010010\dots), (c_n)_0^\infty=(100100\dots)$. Note
first that $a_0$ cannot be equal to 1 nor can  $b_0$, because
$\la<\la_0$ implies $1>\frac{\la}{1-\la^3}>\frac{\la^2}{1-\la^3}$.
Hence $a_0=b_0=0,c_0=1$. Similarly, $a_1=0$ for the same reason as
$b_0$, and $c_1=0$ as well, because $\la>\frac{\la^3}{1-\la^3}$.
Thus, $b_1=1$. Finally, $c_2$ and $b_2$ must be equal to 0, whilst
$a_2=1$.

Thus, we have
\begin{align*}
a_3+a_4\la+a_5\la^2+\dots&=\frac{\la^2}{1-\la^3},\\
b_3+b_4\la+b_5\la^2+\dots&=\frac{\la}{1-\la^3},\\
c_3+c_4\la+c_5\la^2+\dots&=\frac{1}{1-\la^3},
\end{align*}
and we can continue the process {\em ad infinitum}. Therefore, each
$a_k,b_k$ and $c_k$ is uniquely determined from the system of
equations~(\ref{eq:triple}), whence $\bm\pi(\la)\in\U_\la$.

\medskip\noindent (2) Suppose $\la>\la_0$.
Following the argument of the proof of Theorem~\ref{thm:multi-ejk},
we introduce the sets
\[
\De_i:=\De\setminus\bigcup_{j\neq i}f_j(\De)
\]
(three rhombi). Thus, if $\x\in\U_\la$, then necessarily
$\x\in\bigcup_i\De_i$. Fix $i\in\{0,1,2\}$; again, since the shift
map on $\De_i$, i.e., $f_i^{-1}$, is expanding on
$\De_i\setminus\{\p_i\}$, eventually $f_i^{-n}(\x)\in\De_i$ and
$f_i^{-n-1}(\x)\notin\De_i$ for some $n\ge0$, for any
$\x\in\De_i\setminus\{\p_i\}$.

If $f_i^{-n-1}(\x)\notin\bigcup_{j\neq i}\De_j$, then
$\x\notin\U_\la$. Put
\[
\Gamma_i=f_i^{-1}(\De_j)\cap\De_i,\quad j\neq i.
\]
In view of the symmetry, the choice of $j\neq i$ is unimportant --
see Fig~\ref{fig:delta}. Thus, $\U_\la\neq\emptyset$ implies
$\bigcup_i\Gamma_i\neq\emptyset$.

Note that the $\Gamma_i$ are equal for $i=0,1,2$, whence
$\bigcup_i\Gamma_i\neq\emptyset\Leftrightarrow\Gamma_0\neq\emptyset$.
The latter is in fact equivalent to $\la<1/\sqrt2$. Indeed, we have
in barycentric coordinates,
\[
\Gamma_0=\{x<(1-\la)/\la,\,y<1-\la,\,z<1-\la\}.
\]
An open triangle $\{x<a,y<b,z<c\}$ is nondegenerate if and only if
$a+b+c>1$. Hence $(1-\la)/\la+2(1-\la)>1$, which is equivalent to
$\la<1/\sqrt2$. Thus, $\la\in(\la_0,1/\sqrt2)$.

Assume $\x\in\Gamma_0\cap\U_\la$; then $f_0^{-1}(\x)$ has to
intersect $\bigcup_i\De_i$. It is easy to check that
$f_0^{-1}(\Gamma_0)\cap\De_0=\emptyset$, whence, in view of the
symmetry, $f_0^{-1}(\Gamma_0)\cap\De_1\neq\emptyset$. We have

\begin{figure}[t]
    \centering 
\centerline{
\begin{picture}(280,250)
\thicklines \path(0,0)(140,242.2)(280,0)(0,0)
\path(45.5,78.8)(136.5,78.8)(91,0)
\path(189,0)(143.5,78.8)(234.5,78.8)
\path(94.5,163.68)(140,84.87)(185.5,163.68) \thinlines
\dottedline(67.4,116.74)(202.8,116.74)(134.81,0)
\dottedline(77.79,116.74)(145.19,0)
\dottedline(202.8,116.74)(212.6,116.74)
\dottedline(119,45)(160.5,45) \dottedline(202.8,116.74)(207.6,125.4)
\dottedline(77.79,116.74)(72.8,125.4)
\put(64,42){{\Large $\De_1$}} \put(192,42){{\Large $\De_2$}}
\put(132,156){{\Large $\De_0$}}
\put(135,105){$\Gamma_0$}\put(112,64){$\Gamma_1$}
\put(156,64){$\Gamma_2$}
\put(-15,-4){$\p_1$}
\put(283,-4){$\p_2$}
\put(135,248){$\p_0$}
\put(125,36){{\tiny $f_0^{-1}(\Gamma_0)$}}
\put(137,-8){$x$}\put(213,121){$y$}\put(60,121){$z$}
\put(122,48){{\tiny $\mathcal M$}} \thicklines
\put(119.5,45){\circle*{1}}
\end{picture}}
\caption{Triangular case:
$f_0^{-1}(\Gamma_0)\cap\left(\bigcup_i\De_i\right)=\emptyset$.}
    \label{fig:delta}
\end{figure}
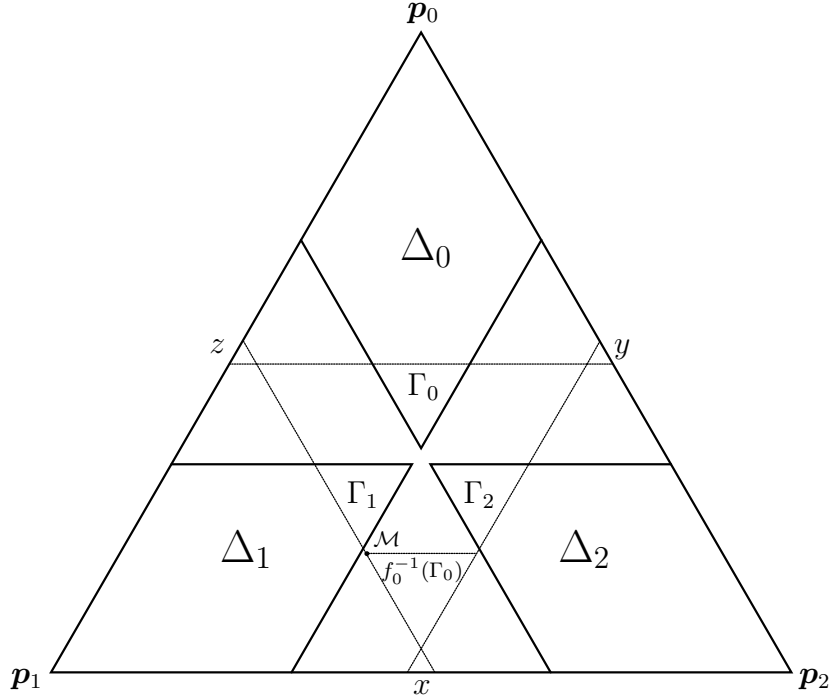

\begin{align*}
f_0^{-1}(\Gamma_0)&=\left\{x<\left(\frac{1-\la}{\la}\right)^2,\,
y<\frac{1-\la}{\la},\,z<\frac{1-\la}{\la}\right\},\\
\De_1&=\{x<1-\la,\,y<1-\la\},
\end{align*}
and we claim that actually, $f_0^{-1}(\Gamma_0)$ lies strictly on
the right of $\De_1$ -- see Fig~\ref{fig:delta}.

Indeed, the coordinates of the point $\mathcal M$, the top left
corner of the triangle $f_0^{-1}(\Gamma_0)$, are as follows:
$\mathcal M\left(\left(\frac{1-\la}{\la}\right)^2,\frac{1-\la}{\la},
\frac{-1+\la+\la^2}{\la^2}\right)$, and we observe that the
inequality $\frac{-1+\la+\la^2}{\la^2}>1-\la$ is equivalent to
$\la>\la_0$.

Thus, each $\x\in\De\setminus\{\p_0,\p_1,\p_2\}$ has at least two
addresses of the form $\x\sim(i_1,\dots,i_n,i_{n+1},\dots)$ and
$\x\sim(i_1,\dots,i_n,j_{n+1},\dots)$ with $j_{n+1}\neq i_{n+1}$ and
a ``compulsory'' prefix $(i_1,\dots,i_n)$, which may be empty (see
Fig.~\ref{fig:branching}). Hence $\x$ has $2^{\aleph_0}$ different
addresses.
\end{proof}

\begin{rmk}One can easily obtain from the proof of the previous
theorem that for $\U_{\la_0}=\emptyset$ as well, but in fact, the
point $\bm\pi(\la_0)=(\la_0^4,\la_0^3,\la_0^2)$ has only $\aleph_0$
different addresses. Thus, $\la_0$ is indeed the full analogue of
the golden ratio for the triangular model. We leave the details as
an exercise for the reader.
\end{rmk}

\begin{rmk}Note that if, like in the proof of
Theorem~\ref{thm:multi-ejk}, $T_\la$ denotes the inverse of
$\{f_0,f_1,f_2\}$ (well defined on $\bigcup_i\De_i$), then
$\bm\pi(\la)$ given by (\ref{eq:pi}) is a period~3 point for
$T_\la$, i.e., $T_\la^3\bm\pi=\bm\pi$. Another 3-cycle is generated
by $\bm\pi'(\la)=\left(\frac{1}{1+\la+\la^2},
\frac{\la}{1+\la+\la^2}, \frac{\la^2}{1+\la+\la^2}\right)$, and we
conjecture that if $\frac23\le\la<\la_0$, then $\bigcup_{i=0}^2
\Gamma_i\cap\U_\la$ consists of just these 6 points. This would
imply that $\U_\la$ is countable for this range of parameters.

A full ``triangular'' analogue of Theorem~\ref{thm:gs} is yet to be
determined. In particular, what is the analogue of the
Komornik-Loreti constant for the triangular case?

Note that for $\la=g=(\sqrt5-1)/2$ the set $\U_\la$ is a continuum
naturally isomorphic to the space of one-sided 0-1 sequences, and
its Hausdorff dimension is $-\log2/\log g$ -- see
\cite[Theorem~6.4]{BMS} and Fig.~5 therein.
\end{rmk}

\medskip\noindent{\bf Acknowledgment.} The author is indebted to
F.~Petrov for his generous help with the proof of
Proposition~\ref{prop:dd} and especially to the anonymous referee
for many useful remarks and suggestions.

\end{document}